\documentstyle[12pt,amssymb]{article}
\begin{document}

\setlength{\parskip}{3pt plus 5pt minus 0pt}

\vskip 1cm
\title{Amenability, Bilipschitz Equivalence, and the Von Neumann Conjecture}
\author{Kevin Whyte, University of Chicago}
\maketitle

\begin{abstract}
  We determine when a quasi-isometry between discrete spaces is at bounded distance from a bilipschitz map.  From this we prove a geometric version of the Von Neumann conjecture on amenability.  We also get some examples in geometric groups theory which show that the sign of the Euler characteristic is not a coarse invariant.  Finally we get some general results on uniformly we finite homology which we will apply to manifolds in a later paper.
\end{abstract}

  In his fundamental paper on amenability, Von Neumann conjectured that if $\Gamma$ is a finitely generated, non-amenable group, then $\Gamma$ has a subgroup which is free on two generators.  While this conjecture proved to be false in general [O], it is true for many classes of groups.  Given the recent interest in studying groups via geometric methods, it seems natural to ask for a geometric version of the conjecture.  If $\Gamma$ does contain a free subgroup then its cosets will partition $\Gamma$ into copies of the free group.  The free group is, geometrically, the regular 4-valent tree.  Our ``geometric Von Neumann conjecture'', is:

 {\bf Theorem 1} {\em If $X$ is a uniformly discrete space of bounded geometry (in particular, for $X$ a finitely generated group or a net in a leaf of a foliation of a compact manifold), $X$ is non-amenable iff it admits a partition with pieces uniformly bilipschitz equivalent to the regular 4-valent tree.}

  The proof depends on constructing a bilipschitz eqivalence of $X \times \{0,1\}$ and $X$ near the projection map.  Our main technical result is a general method for constructing bilipschitz maps near a quasi-isometry.  In particular, we have:

 {\bf Theorem 2} {\em If $X$ and $Y$ are uniformly discrete, non-amenable spaces, $f:X \to Y$ a quasi-isometry, then there is a bilipschitz equivalence between $X$ and $ Y$ at finite distance from $f$.}
      
  Theorem 2 shows, in particular, that $X$ and $Y$ are bilipschitz equivalent.  This gives a generalization of the results in [P], which produces bilipschitz maps between free groups (a question raised by Gromov in [G]).  From this we get:

{\bf Corollary} {\em If $\Gamma_1$ and $\Gamma_2$ are quasi-isometric and non-amenable then, for any $G$, so are $\Gamma_1 * G$ and $\Gamma_2 * G$}

 (A bilipschitz map $\Gamma_1 \to \Gamma_2$ induces one from $\Gamma_1 * G \to \Gamma_2 * G$.)

  This easily gives new examples of groups which are quasi-isometric but do not contain a common subgroup of finite index.  Perhaps more interestingly, one sees that although the sign of the euler characteristic is invariant under passing to subgroups of finite index, it is not invariant under quasi-isometries (this question was also raised by Gromov in [G]).

  If the spaces are amenable, it is no longer true that all quasi-isometries are close to bilipschitz equivalences; for example multiplication by 2 from $\Bbb Z$ to $\Bbb Z$ is clearly not at bounded distance from any surjective map.  To describe the obstruction to finding a nearby bilipschitz map, we need ``uniformly finite homology'' introduced in [BW1] (see section 1 for a definition).  We can completely characterize when it is possible to find a bilipschitz map near a quasi-isometry:

{\bf Theorem 3}  {\em If $X,Y$ uniformly discrete,  $f:X \to Y$ a quasi-isometry, then there is a bilipschitz equivalence between $X$ and $ Y$ at finite distance from $f$ iff $f_*([X]) = [Y]$ in $H_0^{uf}(Y)$.}

  Here $[Z]$ is the class represented by the 0-chain which is sum of all the points in $Z$.  

  Theorem 3 says that the only obstruction is the same sort of density issue that shows up for multiplication by two from $\Bbb Z$ to $\Bbb Z$.  It turns out that if we fix $Y$ but allow $X$ to vary, these obstructions generate $H_0^{uf}(Y)$.  As a consequence of the proof of Theorem 3, we have: 

{\bf Theorem 4} {\em For $c \in C_0^{uf}(X)$, we have $c=0$ in $H_0^{uf}$ iff $$\exists C,r\ s.t\ \forall S \subset X \ |\Sigma_S \ c| \leq C|\partial_r S|$$}

  Where $\partial_r S$ is the set of points outside of $S$ which are within $r$ of S.

  As a corollary of this we get new proofs of many of the results in [BW1], as well as some new facts about uniformly finite homology.  Theorem 4 was asserted without proof in [BW1].  According to Weinberger, their method of proof is completely different.  

  Finally, in a future paper, we will combine the essentially combinatorial results from this paper with some analytic techniques to get new results concerning complete manifolds of bounded geometry.
\\\\

{\bf 1. Definitions and Notation}\\

{\bf Definition} {\em A space $X$ is {\em uniformly discrete} if $\exists r>0$ such that $d(x_1,x_2)<r\ \Rightarrow x_1=x_2$.}

{\bf Definition} {\em A uniformly discrete space $X$ has {\em bounded geometry} if $\forall r$ there is $N_r$ such that for any $x \in X$ $|B_r(x)| \leq N_r$.}

 Unless otherwise stated, all our spaces will be assumed uniformly discrete with bounded geomtery.

{\bf Definition} {\em $f:X\to Y$ is {\em coarse lipschitz} if $\exists A,B$ such that $d(f(x_1),f(x_2)) \leq A d(x_1,x_2) + B $.}

{\bf Definition} {\em $f:X\to Y$ is {\em effectively proper} if $\forall r \ \exists s$ such that $d(f(x_1),f(x_2)) \leq r \ \Rightarrow d(x_1,x_2) \leq s$.}

{\bf Definition} {\em $f:X \to Y$ is a {\em quasi-isometry} (or {\em coarse equivalence}) if $f$ is coarse lipschitz and $\exists g:Y \to X$ coarse lipschitz with $fg$ and $gf$ at bounded distance from the identity maps.} 

 Notice that since the spaces are uniformly discrete, we have that a map is lipschitz if it is coarse lipschitz and injective and likewise a bilipschitz equivalence is just a bijective coarse equivalence.

For $A \subset Z$ we define:

 $$N_r(A) = \{ z \in Z\ s.t.\ d(z,A) \leq r\}$$ and,
 $$\partial_r(A) = N_r(A) - A$$

{\bf Definition} {\em $Z$ is {\em amenable} iff\\ $$\exists \{S_i \subset Z\}\ s.t.\ \forall r>0  \lim_{i\to \infty}\frac {|\partial_r(S_i)|}{|S_i|} =0$$} \\Such a sequence is called a regular sequence.

  The following reformulation will be useful:

{\bf Lemma} The following are equivalent:
\begin{enumerate}
\item $Z$ is non-amenable.
\item $\forall 0<C<1 \exists r\ s.t.\ \forall A\subset Z\ \ |A| \leq C|N_r(A)|$
\item $\exists 0<C<1 \exists r\ s.t.\ \forall A\subset Z\ \ |A| \leq C|N_r(A)|$
\item $\forall C>0 \exists r\ s.t.\ \forall A\subset Z\ \ |A| \leq C|\partial_r(A)|$
\item $\exists C>0 \exists r\ s.t.\ \forall A\subset Z\ \ |A| \leq C|\partial_r(A)|$
\end{enumerate}

{\bf Pf/} \\It is immediate from the definition that all the conditions imply non-amenability.  The last two conditions are just rephrasings of the second and third.  $(2) \Rightarrow (3)$ is trivial, so the two things to prove are $(3) \Rightarrow (2)$ and $(1) \Rightarrow (5)$.

{\bf $(3) \Rightarrow (2)$}\\
$$N_s(N_t(A)) \subset N_{s+t} (A)$$
 so if (3) holds for $C$ and $r$ then:
$$|N_{kr}(A)| \geq |N_r(N_{(k-1)r}(A))| \geq {1/C} | N_{(k-1)r}(A)|$$\\
so,\\
$$|A| \leq C^k |N_{kr}(A)|$$

{\bf $(1) \Rightarrow (5)$}

If (5) does not hold, then there is a sequence $\{S_i\}$ of sets in $Z$ for which the condition fails with $C=i$ and $R=\frac{1}{i}$.  This is then a regular sequence, which contradicts (1).
$\Box$

  To $Z$ a uniformly discrete space we associate a sequence of simplicial complexes $R_r(Z)$.  $R_r(Z)$ has as vertices $Z$, and $(z_0,...,z_k)$ are the vertices of a $k$-simplex iff they are pairwise within $r$ (equivalently, $R_r(Z)$ is the nerve of the covering by $r$-balls.).  We have an inclusion $R_r \to R_{r'}$ for $r' \geq r$.  While the $R_r$ are not individually coarse invariants, the ``pro-system'' is, in particular (and all we need), any sort of limit object or invariant will be coarse.  

{\bf Definition} {\em A chain in $R_r(Z)$ is {\em uniformly locally finite} iff the sum of the norms of the coefficients of simplices in any ball is bounded by some constant depending only on the radius.}

  The group of such chains will be written $C_*^{uf}(Z)$, and the corresponding homology groups are the desired $H_*^{uf}(Z)$.  Unless otherwise mentioned, the coefficients are in $\Bbb Z$, although the definition is sensible for anything with norm (in particular for $\Bbb R$, see section 3).
 
  If $Z$ has bounded geometry, then the sum of all $z \in Z$ is a valid 0-cycle.   We will denote this cycle (and, abusively, the class it represents) by $[Z]$. \\\\

{\bf 2. Finding bilipschitz maps}

{\bf Theorem A} {\em $f:X \to Y$, a coarse equivalence, is at bounded distance from a bilipschitz equivalence iff $f_* ([X])=[Y]$ in $H_0^{uf}(Y)$}

{\bf Pf/}\\
$\Rightarrow$ is obvious.\\
$\Leftarrow$
  Since $X$ and $Y$ are uniformly discrete, a coarse equivalence which is a bijection is a bilipschitz equivalence.  Any map at bounded distance from $f$ is a coarse equivalence, so we need only find a bijection within bounded distance of $f$.  We do this in two steps.  First we find an injective map at near $f$, and one near its ``coarse inverse'' (equivalently, we find a surjection near $f$).  Then we use a Schroeder-Bernstein theorem to produce the bijection.

  The first is handled by the following:

{\bf Lemma} {\em If $f:X \to Y$ is coarse lipschitz, with $f_*([X])=[Y]$ in $H_0^{uf}(Y)$, then there is an injection near $f$.}\\
  
{\bf Pf/} \\
  As before, $\Rightarrow$ is clear.\\  
  $\Leftarrow$   Fix $r>0$.  If we are to move $f$ by no more than $r$, the set of possible values at $x \in X$ is $B_r(f(x))$.  we want to know whether we can choose a value at each $x$ without repetition.  This is exactly the sort of question addressed by Hall's ``Marriage Lemma'':

 {\bf Lemma (Marriage Lemma)} {\em If we have $g:A \to Fin(B)$ (the finite subsets of $B$), then there's an injective map $\phi:A \to B$ with $\phi(a) \in g(a) \forall a \in A$ iff $\forall S \subset A$ we have that}
   {$\displaystyle |S| \leq |\cup_{s \in S} g(s)|$}

 In the case of the lemma, we need $\forall S \subset X$ $|S| \leq |N_r(f(S))|$.  The LHS is increased if we replace $S$ by $f^{-1}(f(S))$.  Therefore we need, 
$$(*)\forall S \subset Y\ \ |f^{-1}(S)| \leq |N_r(S)| $$ 
We have
 $$|f^{-1}(S)| = |\Sigma_{S} f_*([X])|$$
where $\Sigma_S c$, for $c \in C_0^{uf}$, is just the sum of the coefficients of vertices in $S$  

 So $(*)$ becomes 

 $$\Sigma_S f_*([X]) \leq |N_r(S)|$$
equivalently,
 $$\Sigma_S (f_*([X]) - [Y]) \leq |\partial_r(S)|$$

 Let $c=f_*([X]) - [Y]$.  We have that $c \geq -1$ everywhere, and since $f$ is coarsely onto, $c \geq 0$ on a set of density $\geq \lambda$ for some $\lambda > 0$.  By assumption, $c=\partial b$ for some $b$, with the coefficients of $b$ bounded by $M$ for some $M$, and a maximum edge length of $l$.  Consider $\partial_r S$.  The edges of $b$ make this an oriented graph, which has a flow defined by the coefficients of $b$.  This flow has sinks in $S$ with total flow $\Sigma_S c$.  The sources are of two kinds: those at distance less than $r$ from $S$ and those at distance $r$ from $S$.  The former produce at most 1 each, while those at the boundary can produce up to $M$.  If the interior sources produce more than $\lambda$ of the total flow into $S$ then, as their density is bounded by $\lambda$, there are at least $\Sigma_S c$ vertices in $\partial_r S$ as desired.  If the interior sources produce less than $\lambda$ of the total, then at least $(1-\lambda)\Sigma_S c$ must flow in from a distance of $r$. Since the maximum edge length is $l$, and no edge can have more than a flow of $M$, this requires at least $\frac{(1-\lambda)(\Sigma_S c)r}{Ml}$ edges.  As the cycle $b$ is uniformly finite, there is some $v$ so that no vertex has valence more than $v$, thus there must be at least $\frac{(1-\lambda)(\Sigma_S c)r}{Mlv}$ vertices.  Taking $r \geq \frac{Mlv}{(1-\lambda)}$ then gives the required number of vertices. $\Box$   

 Let $g: Y\to X$ be a coarse inverse for $f$.  Then $f_*([X])=[Y]$ gives $(gf)_*([X])=g_*([Y])$, and as $gf$ is at bounded distance from $1_X$, this is $[X]=g_*([Y])$.  So, under the hypotheses of Theorem A, we can produce injections near $f$ and $g$.  Thus we are done with the first step of the argument.

  The standard proof of the Schroeder-Bernstein theorem (see, for example, [H]) actually proves something slightly stronger than is usually stated.

{\bf Theorem}(Schroeder-Bernstein) If $f:A \to B$ and $g:B\to A$ are injective, then there is a bijection from $A \to B$ which agrees everywhere either with $f$ or $g^{-1}$.

 In our case, since $g$ is a coarse inverse for $f$, $g^{-1}$ (where it is defined), is at bounded distance from $f$.  Thus the bijection produced by Schroeder-Bernstein is at bounded distance from $f$.$\Box$

{\bf 3.$H_0^{uf}$ and Amenability}

{\bf Definition} {\em $S_{lip}(X)$ is the set of pairs $(Z,f)$, where $Z$ is uniformly discrete of bounded geometry and $f: Z \to X$ is a coarse equivalence, up to the relation that $(Z_1,f_1)$ is equivalent to $(Z_2,f_2)$ iff there is a bijection $Z_1 \to Z_2$ which commutes (up to bounded distance) with \{$f_i$\}.}

{\bf Definition} {\em A class $c \in H^{uf}_0$ is {\em positive} if there is a chain representing it which is everywhere non-negative, and positive on a C-dense set.}

 {\bf Theorem B} {\em $S_{lip}(X)=H^{uf,+}_0(X)$ (the positive part of $H^{uf}_0$).}

 {\bf Pf/}\\ The map is given by $(Z,f) \mapsto f_*([Z])$.  The map is injective by Theorem A.  To see it is surjective, let $c$ be a positive cycle.  \\Define $$Z_c = \{(x,k) \in X\times{\Bbb N}, k < c_x\}$$
 $(Z_c,\pi_Z) \in S_{lip}(X)$ represents $c$.$\Box$

 Notice that this gives $S_{lip}(X)$ a semi-group structure.  From the proof we can see that the addition is given geometrically by ``parallel disjoint union over $X$''.

  We can apply these methods to give results about $H_0^{uf}$.  Our main result in this direction is the following:

{\bf Theorem C} {If $c$ is a cycle in $C_0^{uf}(Z)$, then $[c]=0$ in $H_0^{uf}(Z)$ iff

 $$\exists C,r\ \forall S \subset Z\ |\Sigma_S c| \leq C |\partial_r S| $$

 This result has also been proven by Block and Weinberger (unpublished) by different methods.  

 {\bf Pf/}

  The condition is clearly necessary.  Since $c$ is uniformly finite, there is some $M$ such that all coefficients of $c$ are less than $M$ in absolute value.  For $n \geq M$, let $X_n$ be the space in $S_{lip}(Z)$ representing $c+n[Z]$ and $Y_n$ the space representing $n [Z]$.  By the last result $c$ bounds iff there is a bilipschitz equivalence between these spaces commuting, up to bounded distance, with the projections to $Z$.  We will build such a map as in the proof of Theorem A.  Fix $r$.  $f:X_n \to Y_n$ commutes within $r$ with the projections iff for all $x,\ f(x)\in \{y:\ d(\pi_X(x),\pi_Y(y)) \leq r\}$.

 Thus, for any set $S$ in $X_n$, the possible range is $\pi_Y^{-1}(N_r(\pi_X(S))$.  To find an injective map the Marriage lemma requires that this be no smaller than $S$.  This is hardest when $S=\pi_X^{-1}(\pi_X(S))$.  We can then write things in terms of $A=\pi_X(S)$, when the condition becomes:

 $$ \forall A \subset Z \ \ |\pi_X^{-1}(A)| \leq |\pi_Y^{-1}(N_r(A))| $$
 In terms of $c$, this says:
 $$ \forall A \subset Z\ \ \Sigma_A c+n[Z] \leq \Sigma_{N_r A} n[Z] $$
 or, 
 $$ \forall A \subset Z\ \ \Sigma_A c \leq \Sigma_{\partial_r A} n[Z] $$
 or
 $$ \forall A \subset Z\ \  \Sigma_A c \leq n|\partial_r A| $$
 So we can build an injection $X_n$ to $Y_n$, for $n$ large enough, iff 
 $$\exists C,r\ s.t.\ \forall S \subset Z\ \ \Sigma_S c \leq C |\partial_r S| $$
 To build the injection the other way, we need
 $$ \forall A \subset Z\ \ \Sigma_A (n[Z]) \leq \Sigma_{N_r A} (c+n[Z]) $$
 equivalently,
 $$ \forall A \subset Z\ \ \Sigma_A(-c)\leq \Sigma_{\partial_r A} (c+n[Z]) $$
 the RHS is at least $(n-M)|\partial_r A|$, and no more than $(n+M)|\partial_r A|$ so the injection can be built this direction, for $n$ large, iff 
 $$\exists C,r\ s.t.\ \forall S \subset Z\ \ \Sigma_S (-c) \leq C |\partial_r S| $$
 The two conditions together give the Theorem.$\Box$  
\\

  This result has applications to characteristic numbers for manifolds with bounded geometry (in particular, as in [BW2], for  manifolds of positive scalar curvature) which we will explore in another paper.

 {\bf Corollary} $H^{uf}_0(Z,{\Bbb Z}) \to H^{uf}_0(Z,{\Bbb R})$ is injective.

 The inequality in Theorem C follows immediately from $c=\partial b$ even if $b$ is allowed real coefficients.\\

 {\bf Corollary} $H^{uf}_0(Z,{\Bbb Z}) \to H^{uf}_0(Z,{\Bbb R})$ is an isomorphism if $Z$ is infinite and coarsely connected.

 (Coarse connectivity means that $\exists r$ such that any two points of $x,y \in Z$ can be joined by a chain $x=z_0,z_1,z_2, \dots , z_n=y$ with $d(z_i,z_{i+1})<r$.)  

 {\bf Pf/}

   We just saw injectivity.  For surjectivity:

   From the sequence $0\to{\Bbb Z} \to {\Bbb R} \to S^1\to 0$ we get a long exact sequence in $H^{uf}_*$ which ends with:

   $\dots \to H^{uf}_0(Z,{\Bbb Z}) \to H^{uf}_0(Z,{\Bbb R}) \to H^{uf}_0(Z,S^1) \to 0$

  Since $S^1$ is bounded, locally finite $\Rightarrow$ uniformly locally finite.  Thus,

   $H^{uf}_0(Z,S^1) = \lim_{r \to \infty} H_0^{lf}(R_r(Z),S^1)$

 We have:

   $\dots \to H^{lf}_0(R_r(Z),{\Bbb Z}) \to H^{lf}_0(R_r(Z),{\Bbb R}) \to H^{lf}_0(R_r(Z),S^1) \to 0$

 If $Z$ is coarsely connected and infinite, $H^{lf}_0(R_r(Z),{\Bbb R})=0$ for $r$ large enough.  Therefore $H^{uf}_0(Z,S^1)=0$, proving surjectivity. $\Box$

 So, if $H_0^{uf}(Z,{\Bbb Z})$ is non-zero, it's uncountable and even has an ${\Bbb R}$ vector space structure.

  The inequality in Theorem C is true for any bounded chain if we have that $\forall C>0 \exists r \forall S \subset Z\ \  |S| \leq C|\partial_r(S)|$.  This is exactly one of the equivalent conditions for non-amenability given in the first section.  On the other hand, for a positive cycle to vanish, we must have this condition, therefore:

{\bf Theorem D} {\em The following are equivalent:
\begin{enumerate}
\item $X$ is non-amenable 
\item $H_0^{uf}(X,{\Bbb Z}) = 0$
\item $H_0^{uf}(X,{\Bbb R}) = 0$
\item $[X]=0$ in $H_0^{uf}(X,{\Bbb Z}) = 0$
\item $0 \in H_0^{uf,+}(X,{\Bbb Z}) $
\item $S_{lip}(X)= \{pt\}$
\end{enumerate}}
 All but the last are proven analytically in [BW2].

The last condition is a characterization of ``strongly lipschitz rigid'' spaces, where strongly refers to the fact that get rigidity for maps, not just spaces.  ``Weak lipschitz rigidity'' for $X$ should mean that any space coarsely equivalent to $X$ is bilipschitz to $X$.  Any non-amenable space is weakly rigid, but the converse is not true: ${\Bbb Z}$ is amenable and weakly rigid.

{\bf Corollary} {\em $X$ is ``weakly rigid'' iff the group of quasi-isometries of $X$ is transitive on $H_0^{uf,+}(X)$}

  This is unsatisfactory in that the group of self quasi-isometries is rarely a tractable object.  It is not hard to produce examples of non-weakly rigid spaces by building spaces with no non-trivial quasi-isometries (for example, $\{n^2\} \subset {\Bbb Z}$), but symmetric examples seem hard.  It has only recently been determined that ${\Bbb Z}^2$ is not weakly rigid, is open [BK].

{\bf 4.  The Geometric Von Neumann Conjecture}

{\bf Theorem} (Geometric Von Neumann Conjecture) {\em $Z$ is non-amenable iff $Z$ has a partition whose pieces are lipschitz embedded copies of the tri-valent tree (with uniform lipschitz constant)}

  We note that we can just as well have any tree quasi-isometric to the tri-valent tree as components.  This follows from Theorem A since we can use a bijection between the trees to change one to the other.  In particular, we can take the uniformly 4-valent tree for the components.  The Von Neumann conjecture (which is false) is that any finitely generated, discrete, non-amenable group contains a free subgroup of rank two.  The coset space of this subgroup would then give a partition of the type demanded by the theorem.

{\bf Pf/}

  ``if'' is clear.

  For ``only if'', we construct a partition with components (uniformly) quasi-isometric to the tri-valent tree.  As in the above discussion, this will prove the theorem.  

{\bf Lemma} {\em If $Z$ is non-amenable then there is a bilipschitz map $f$ from $Z \times \{0,1\}$ to $Z$ near the projection map but never equal to it.}

 {\bf Pf/} (of lemma)

 Except for the last bit, this is immediate from the earlier theorems.  The set of possible values for $f(z,\epsilon)$, for $f$ satisfying the conditions in the lemma and $(z,\epsilon) \in Z \times \{0,1\}$ is $B_r(z) - \{z\}$.  So, for any $S$ in $Z$, the set of possible values for $f(S \times \{0,1\})$ is at least $\partial_r(S)$.  By non-amenability, $|\partial_r(S)|$ can be made larger than $2|S|$ by taking $r$ large enough.  Thus, as in Theorem A, we can produce an injection with the desired properties.  Schroeder-Berstein does the rest.   

 Given such a bilipschitz map $f$, take $Z$ to be the set of vertices of a graph whose edges are $(z,f(z,0))$ and $(z,f(z,1))$.  Thus every vertex has valence three (the obvious two ``outgoing'' edges, and the ``incoming'' $(f^{-1}(z),z)$).  Call this graph $\Gamma$.

{\bf Claim} Each connected component of $\Gamma$ has at most one loop.

  Consider a path in $\Gamma$.  Since each vertex has only one incoming edge, any path starting with an oriented edge must be oriented.  Thus any path is either oriented, or there is some vertex with both of its outgoing edges on the path.  In this case the path is the union of two oriented paths leaving that point.  Such a path cannot form a loop as the point where the two parts of the path re-joined would need two incoming edges.  Thus, if there is a loop, it is oriented.  If there were two loops in a component, we would have that the path connecting them would point away from both loops (since the incoming edges for the endpoints of the path are already used in the loops).  This would force a point in the middle of the path to have two incoming edges, which is impossible.

  Cutting at most one edge in each component, we get a decomposition into trees.  These trees are tri-valent, except at at most two vertices, hence uniformly bilipschitz to the tri-valent tree. $\Box$

{\bf 5.Free Products}

   Given two finitely generate discrete groups, one can form the free product.  While this is a natural thing to do algebraically (and homotopically on the level of $K(\Gamma,1)$s), it is poorly behaved coarsely.  For example, while ${\Bbb Z}/2{\Bbb Z}$ is coarsely equivalent to \{1\}, it is definitely not the case that ${\Bbb Z}/2{\Bbb Z}*{\Bbb Z}/2{\Bbb Z}$ is coarsely equivalent to  $\{1\}*\{1\}$ (or to ${\Bbb Z}/2{\Bbb Z} * \{1\}$).

 {\bf Theorem} {\em If $\Gamma_1$ and $\Gamma_2$ are coarsely equivalent finitely generated non-amenable groups, then for any $G$, $\Gamma_1 * G$ and $\Gamma_2 * G$ are bilipschitz equivalent.}

{\bf Pf/}

  By the earlier theorems, we can find a bilipschitz map (WLOG taking 1 to 1) from $\Gamma_1$ and $\Gamma_2$.  This gives a bilipschitz equivalence between the corresponding free products by applying the map term by term to reduced words. $\Box$

  Note that this equivalence between the free products is not well determined (even up to bounded distance) - it depends on the choice of bilipschitz map near the original quasi-isometry.  

{\bf Question} {\em Are there finitely generated torsion free groups which are quasi-isometric, but do not remain so after taking free products with a fixed group?}

  Since all surface groups (of higher genus) are quasi-isometric and non-amenable, they are bilipschitz.  Their free products are then also bilipschitz.  The free product $\Sigma_2 * \Sigma_3$ contains $\Sigma_2 * \Sigma_2 * \Sigma_5$ as a subgroup of index 2.  This shows that the free product of two surface groups is also quasi-isometric (and hence bilipschitz) to the free product of three such groups.  Thus, all free products of at least two surface groups are quasi-isometric.  These groups are almost never commensurable; for example, it is an easy consequence of the Kurosh subgroup theorem that $\Sigma_a * \Sigma_a$ and $\Sigma_b * \Sigma_b$ cannot have a common subgroup of finite index.

  Consider the groups $\Gamma=\Sigma_2 \times \Sigma_2$ and $\Delta=\Sigma_3$.  Let $\Gamma_k$ be a subgroup of $\Gamma$ of index $k$, and $\Delta_j$ a subgroup of $\Delta$ of index $j$.  By the applications of the theorem, all free products of the form $\Gamma_k * \Delta_j$ are quasi-isoemtric.   The euler charactistic of $\Gamma_k * \Delta_j$ is $4(k-j) - 1$.  Thus, by varying $j$ and $k$ this euler characteristic can be either positive or negative.  This answers negatively a question of Gromov in [G] as to whether this sign is a quasi-isometry invariant.

{\bf References}
\newline
$[BK]$ Burago and Kleiner, {\em Separated nets in Euclidean space and Jacobians of bilipschitz maps} pre-print.\\
$[BW1]$ J.Block and S.Weinberger, {\em Large Scale Homology Theories and Geometry} in AMS/IP Studies in Advanced Mathematics, Vol 2 (1997).\\
$[BW2]$ J.Block and S.Weinberger, {\em Aperiodic Tilings, Positive scalar curvature, and Amenability }, JAMS 5 (1992).\\
$[G]$ M. Gromov {\em Asymptotic Invariants of Infinite Groups}, in Geometric Group Theory, (Nibling and Roller,eds), LMS lecture notes v.182 (1993).\\
$[H]$ P. Halmos, {\em Naive Set Theory}, Springer-Verlag. \\
$[O]$ Olshanskii, {\em On the question of the existence of an invariant mean on a group}, Uzpekhi Mat. Nauk. 35 (1980).\\
$[P]$ P. Papasoglu, {\em Homogeneous Trees are Bi-Lipschitz Equivalent}, Geom. Dedicada, 54 (1995).\\  

\end{document}